\documentclass[12pt]{amsart}
\usepackage{amsmath,amscd,amssymb,amsfonts,graphics}
\newcommand{\h}{\hbox}
\newcommand{\q}{\quad}
\newcommand{\nin}{\noindent}
\newcommand{\bs}{\par\bigskip}
\newcommand{\ms}{\par\medskip}
\newcommand{\sk}{\par\smallskip}
\newcommand{\bsn}{\par\bigskip\noindent}
\newcommand{\msn}{\par\medskip\noindent}

\newcommand{\ges}{\geqslant}
\newcommand{\les}{\leqslant}
\newcommand{\1}{\hskip1pt}
\newcommand{\mcap}{\hbox{$\bigcap$}}
\newcommand{\mcup}{\hbox{$\bigcup$}}
\newcommand{\msqcup}{\hbox{$\bigsqcup$}}
\newcommand{\msum}{\hbox{$\sum$}}
\newcommand{\mopl}{\hbox{$\bigoplus$}}

\newcommand{\I}{{\mathcal I}}
\newcommand{\K}{{\mathcal K}}

\newcommand{\Hc}{{\mathcal H}}
\newcommand{\M}{{\mathcal M}}
\newcommand{\Nc}{{\mathcal N}}
\newcommand{\OO}{{\mathcal O}}

\newcommand{\X}{{\mathcal X}}

\newcommand{\C}{{\mathbb C}}
\newcommand{\DD}{{\mathbb D}}
\newcommand{\PP}{{\mathbb P}}
\newcommand{\Q}{{\mathbb Q}}
\newcommand{\Lb}{{\mathbb L}}
\newcommand{\N}{{\mathbb N}}
\newcommand{\R}{{\mathbb R}}
\newcommand{\RR}{{\mathbf R}}
\newcommand{\TT}{{\mathbb T}}
\newcommand{\Z}{{\mathbb Z}}

\newcommand{\Gr}{{\rm Gr}}
\newcommand{\al}{\alpha}

\newcommand{\Ga}{\Gamma}

\newcommand{\De}{\Delta}

\newcommand{\si}{\sigma}
\newcommand{\Si}{\Sigma}
\newcommand{\ep}{\varepsilon}
\newcommand{\Xo}{{}\,\overline{\!X}{}}
\newcommand{\Zo}{{}\,\overline{\!Z}{}}
\newcommand{\Om}{\Omega}
\newcommand{\om}{\omega}
\newcommand{\dd}{\partial}
\newcommand{\ddd}{{\rm d}}

\newcommand{\Zf}{Z_{\!\1f}}
\newcommand{\Zfs}{Z_{\!\1f,\si}}
\newcommand{\Zof}{\Zo_{\!f}}
\newcommand{\Zofs}{\Zo_{\!f,\si}}
\newcommand{\MHM}{{\rm MHM}}
\newcommand{\MHS}{{\rm MHS}}

\newcommand{\pl}{\,{+}\,}
\newcommand{\mi}{\,{-}\,}
\newcommand{\eq}{\,{=}\,}
\newcommand{\bl}{\bigl}
\newcommand{\br}{\bigr}

\newcommand{\ssb}{\raise.15ex\h{${\scriptscriptstyle\bullet}$}}
\newcommand{\ssc}{\,\raise.15ex\h{${\scriptstyle\circ}$}\,}

\newcommand{\into}{\hookrightarrow}
\newcommand{\simto}{\,\,\rlap{\hskip1.5mm\raise1.4mm\hbox{$\sim$}}\hbox{$\longrightarrow$}\,\,}
\newcommand{\nsset}{\rlap{\raise.5mm\h{$\subset$}}\raise-1.5mm\h{$\rlap{\raise.5mm\h{$\scriptscriptstyle\,\,/$}}-$}}
\begin{document}
\title[Intersection complexes]
{Intersection complexes of toric varieties\\and mixed Hodge modules}
\author[M. Saito]{Morihiko Saito}
\address{M. Saito : RIMS Kyoto University, Kyoto 606-8502 Japan}
\email{msaito@kurims.kyoto-u.ac.jp}
\begin{abstract} We prove the structure theorem of the intersection complexes of toric varieties in the category of mixed Hodge modules. This theorem is due to Bernstein, Khovanskii and MacPherson for the underlying complexes with rational coefficients. As a corollary the Euler characteristic Hodge numbers of non-degenerate toric hypersurface can be determined by the Euler characteristic subtotal Hodge numbers together with combinatorial data of Newton polyhedra. This is used implicitly in an explicit formula by Batyrev--Borisov. Note that a formula for the Euler characteristic subtotal Hodge numbers in terms of Newton polyhedra has been given by Danilov--Khovanskii. The structure theorem also implies that the graded quotients of the weight filtration on the middle cohomology of the canonical compactification of a non-degenerate toric hypersurface have Hodge level strictly smaller than the general case except for the middle weight. This gives another proof of a formula for the frontier Hodge numbers of non-degenerate toric hypersurfaces due to Danilov--Khovanskii.
\end{abstract}
\maketitle
\centerline{\bf Introduction}
\bsn
For a Laurent polynomial $f$ of $n$ variables, set $Z_{\!f}:=f^{-1}(0)\subset(\C^*)^n$, and $\De:=\Ga(f)$, the Newton polyhedron of $f$, where we assume $n\ges 2$, see (2.1) below. We say that the toric hypersurface $\Zf$ is {\it non-degenerate\1} (or more precisely, {\it stratified-smooth\1}) if the hypersurfaces $f_{\si}^{-1}(0)\subset(\C^*)^n$ is {\it smooth\1} for any face $\si\les\De$ including the case $\si=\De$, where $f_{\si}$ is the partial sum for monomials of $f$ supported on $\si$. This definition is slightly different form \cite{Ko}, where the condition is {\it not\1} imposed for $\De$ itself, see Remark~(2.1) below.
\sk
A non-degenerate toric hypersurface is a natural generalization of smooth hypersurfaces in $\PP^n$, where the Hodge numbers were calculated by Griffiths and Hirzebruch using the Jacobian rings of $f$ and Hirzebruch characteristic classes respectively. It does not seem easy to generalize these directly to non-degenerate toric hypersurfaces, since the canonical compactification $\Zof$ in the toric variety $\X$ defined by the dual fan $\Si$ of $\De$ may have certain {\it singularities\1} coming from the ones of $\X$, see (2.1) below.
\sk
It seems rather easy to calculate the Euler characteristic {\it subtotal\1} Hodge numbers with compact supports
$$\chi^p_c(\Zf):=\msum_j\,(-1)^j\dim_{\C}\Gr_F^pH_c^j(\Zf),$$
in terms of the Newton polyhedron $\De$ as is explained by {Danilov and Khovanski\u\i} \cite{DK}. We can deduce from these the Euler characteristic Hodge numbers with compact supports
$$\chi^{p,q}_c(\Zf):=\msum_j\,(-1)^j\dim_{\C}\Gr_F^p\Gr^W_{p+q}H_c^j(\Zf),$$
using a ``prime cutting" of $\De$, see \h{\it loc.\,cit.} and (2.6) below. Here we call $\chi^p_c(\Zf)$ {\it subtotal,} since $\chi^p_c(\Zf)=\sum_q\chi^{p,q}_c(\Zf)$.
\sk
By the well-known estimates of weights (see (1.4.4) below) and the weak Lefschetz type theorem (or Artin's vanishing theorem \cite{BBD}, see also (1.4.8) below) together with the {\it ampleness\1} and {\it stratified smoothness\1} of $\Zof\subset\X$ (see Remark~(2.4) below), we can get the following equalities :
$$\Gr^W_kH_c^j(\Zf)=\begin{cases}\Gr^W_{k+2}H_c^{j+2}(X_{\De})(1)\q\q\q\q\,\,\,(j\1{>}\1n{-}1),\\ \,0\q\q\q(j\1{<}\1n{-}1\,\,\,\h{or}\,\,\,j\1{=}\1n{-}1,\,k\1{>}\1n{-}1),\end{cases}
\leqno(1)$$
with $X_{\De}\cong(\C^*)^n$. Thus only the $\Gr^W_kH_c^{n-1}(\Zf)$ ($k\,{\les}\,n{-}1$) are non-trivial.
\sk
It is known that $\Gr_{n-1}^WH_c^{n-1}(\Zof)$ is closely related to the {\it intersection complex\1} ${\rm IC}_{\Zof}\Q$, which is isomorphic to the restriction of ${\rm IC}_{\X}\Q[-1]$ to $\Zof$ as a consequence of the {\it stratified smoothness\1} of $\Zof$. We can calculate ${\rm IC}_{\X}\Q$ using a formula due to Bernstein, {Khovanski\u\i} and MacPherson, see \cite{St}, \cite{DL}, \cite{Fi}, etc. More precisely, we have the following.
\msn
{\bf Theorem~1.} {\it Let $\X$ be the toric variety associated with the dual fan $\Si$ of an $n$-dimensional polyhedron $\De\subset\R^n$ defined by a finite number of linear inequalities with rational coefficients.
For any face $\si\les\De$, there are non-canonical isomorphisms of $\Q$-local systems on $X_{\si}$
$$\Hc^{j-n}\bl({\rm IC}_{\X}\Q|_{X_{\si}}\br)\cong\begin{cases}\Q_{X_{\si}}^{m_{\si,j}}(-k)&\h{if}\,\,\,j=2k<c_{\si},\,\,k\in\N,\\
\,0&\h{otherwise.}\end{cases}
\leqno(2)$$
These are liftable to isomorphisms of mixed Hodge modules on $X_{\si}$ {\rm(}up to a shift of complexes, see Remark~{\rm(1.4)\,(iii)} below{\rm)}. Here $X_{\si}$ is a stratum of the stratification $\X=\msqcup_{\si\les\De}\,X_{\si}$ by torus action orbits indexed by faces $\si\les\De$, and $c_{\si}:={\rm codim}\,\si={\rm codim}\,X_{\si}$. Moreover, the ranks of local systems $m_{\si,j}$ are calculated recursively using the combinatorial data of $\De=\Ga(f)$, see $(1.7.9)$ below.}
\ms
The {\it liftability\1} of the isomorphisms in (2) does not seem to be stated explicitly in the literature. Notice a non-trivial difference between $\Hc^{j-n}\bl({\rm IC}_{\X}\Q|_{X_{\si}}\br)$ and $(\Hc^{j-n}{\rm IC}_{\X}\Q)|_{X_{\si}}$ from the view point of mixed Hodge modules. The former is defined by $H^{j-c_{\si}}i_{X_{\si}}^*{\rm IC}_{\X}\Q_h$ (up to a shift of complex, see see Remark~(1.4)\,(iii) below), where $i_{X_{\si}}:X_{\si}\into\X$ is the inclusion, and ${\rm IC}_{\X}\Q_h$ is a pure Hodge module of weight $n$ whose underlying $\Q$-complex is ${\rm IC}_{\X}\Q$. Note that $H^{j-c_{\si}}i_{X_{\si}}^*$ is the {\it composition\1} of the pull-back functor $i_{X_{\si}}^*$ between the bounded derived categories of mixed Hodge modules with the standard cohomological functor $H^{j-c_{\si}}$ of the derived categories. It is rather complicated to define a bounded complex of mixed Hodge modules whose underlying $\Q$-complex is $\Hc^{j-n}{\rm IC}_{\X}\Q$, see \cite[Remark~4.6,\,2]{mhm}.
\sk
It has been found quite recently that Theorem~1 follows easily from the decomposition theorem for a {\it toric blow-up\1} of $\X$ along $X_{\si}$ using the purity and the hard Lefschetz theorem for the intersection cohomology of the exceptional divisor $Y$, where we may assume $\dim\si=0$. This gives a proof of a recursive formula due to Bernstein, {Khovanski\u\i} and MacPherson. Although there is a slight difference between the formulations in \cite{St} and in \cite{DL}, \cite{Fi}, they are essentially the same, see (1.7) below. It is also possible to use the {\it generalized Thom-Gysin sequence\1} for the intersection cohomology of $\X\setminus X_{\si}$ which is a $\C^*$-bundle over $Y$, assuming that $\X$ is an affine toric variety (which is an affine cone of $Y$ with vertex $\si$). Here one can pass to the {\it Grothendieck ring\1} of mixed Hodge structures, since the degree of cohomology is uniquely determined by the weight, and moreover the truncation by degree can be replaced with that by weights, see (1.7) below. (This kind of argument may be needed also in the $\ell$-adic argument.)
\sk
From Theorem~1 we can deduce the following.
\msn
{\bf Corollary~1.} {\it The Euler characteristic Hodge numbers $\chi_c^{p,q}(\Zof)$ {\rm(}and hence $\chi_c^{p,q}(\Zf))$ can be determined from the Euler characteristic subtotal Hodge numbers $\chi_c^p(\Zf)$ using combinatorial data of $\De=\Ga(f)$.}
\ms
Indeed, there are isomorphisms ${\rm IH}^j(\Zof)={\rm IH}^j(\X)$ ($j<n{-}1$) by the weak Lefschetz type theorem, since $\Zof\subset\X$ is an {\it ample\1} divisor and is stratified-smooth. Combining this with the hard Lefschetz theorem, only the middle intersection cohomology ${\rm IH}^{n-1}(\Zof)$ is non-trivial, and its weight is $n{-}1$ by the purity. We then get Corollary~1 by induction on $\dim\X$ using Theorem~1 and the stratification $\Zof=\msqcup_{\si\les\De}\,\Zfs$ with $\Zfs:=\Zof\cap X_{\si}$, since
$${\rm IC}_{\Zof}\Q={\rm IC}_{\X}\Q|_{\Zof}[-1],\q{\rm IC}_{\Zof}\Q|_{\Zf}=\Q_{\Zf}[n{-}1].$$
\sk
There is another way of determining $\chi_c^{p,q}(\Zf)$ replacing the intersection cohomology ${\rm IC}_{\Zof}\Q$ in the above argument with the cohomology of the compactification $\Zof'$ of $\Zf$ in a toric quasi-resolution $\X'$ of $\X$ associated with a ``prime cutting" of $\De$ in (1.3) below, see \cite[5.2]{DK} and (2.6) below.
This algorithm does not seem to be mentioned in remarks related to {\it loc.\,cit.}\ around \cite[Prop. 3.8]{BB}.
It does not seem to give a much simpler formula than the one in {\it loc.\,cit.} Note that the fibers of a quasi-resolution associated with a ``prime cutting" of $\De$ are identified with quasi-resolutions of the projective varieties used by \cite{DL}, \cite{Fi} in the recursive formula for the intersection complexes of toric varieties, see (1.7) below.
\sk
Corollary~1 may be used implicitly in an {\it explicit\1} formula for the $\chi_c^{p,q}(\Zf)$ in \cite{BB}. Here one can {\it embed\1} $\De$ in $\R^{n+1}$ using the inclusion $\R^n=\R^n{\times}\{1\}\subset\R^{n+1}$, and consider the cone $C_{\De}$ of $\De$ in $\R^{n+1}$ so that the ``Gorenstein" property is satisfied {\it automatically\1} (and {\it no condition\1} is imposed on $\De$). The ``degree" is then given by the last coordinate of $\R^{n+1}$, and it corresponds to $k$ of $k\De$. (These do not seem to be mentioned explicitly in {\it loc.\,cit.}) Here one seems to need also the Ehrhart-Macdonald reciprocity \cite{Ma} between the Ehrhart polynomial and the interior Ehrhart polynomial. (The latter seems to follow also from the Hirzebruch-Riemann-Roch theorem \cite{Hi} and Serre duality on a desingularization of the toric variety associated with a smooth subdivision of the fan, since the dualizing sheaf is identified with the ideal sheaf vanishing along the complement of the torus in $\X$, see for instance \cite{DK}.)
\ms
For an algebraic variety $Y$ in general, set
$$h_c^{j,p,q}(Y):=\dim_{\C}\Gr_F^p\Gr^W_{p+q}H_c^j(Y).$$
From the isomorphisms in (2), we can deduce also the following.
\msn
{\bf Corollary~2.} {\it For $k\1{<}\1n{-}1$, the Hodge level of $\Gr^W_kH_c^{n-1}(\Zof)$ is strictly smaller than $k$, that is,}
$$h_c^{n-1,p,q}(\Zof)=0\q\h{\it if}\q p{+}q<n{-}1,\,p\1q=0.
\leqno(3)$$
\ms
Applying this to $\Zofs=\Zof\cap\Xo_{\si}$ for any $\si\les\De$, we can reprove a formula of Danilov and {Khovanski\u\i} \cite{DK} for the {\it frontier} Hodge numbers of a non-degenerate toric hypersurface $\Zf$\,:
$$h_c^{n-1,p,0}(\Zf)=\begin{cases}\msum_{d_{\si}=p+1}\,l^*(\si)&(p>0),\\ \Pi{-}1&(p=0).\raise5mm\h{}\end{cases}
\leqno(4)$$
Here $l^*(\si):=|\si^{\circ}\cap\Z^n|$ with $\si^{\circ}$ the relative interior of $\si$, and $\Pi:=\bl|\mcup_{d_{\si}=1}\,\si\cap M\br|$ with $d_{\si}:=\dim\si$, see (2.5) below.
\sk
This work was partially supported by JSPS Kakenhi 15K04816.
I thank A.~Stapledon for a remark about the algorithm of Danilov--{Khovanski\u\i} \cite{DK}.
\sk
In Section~1 we review some basics of toric varieties, mixed Hodge modules, intersection complexes, and prove Theorem~1. In Section~2 we explain some applications of Theorem~1 to the study of the cohomology of non-degenerate toric hypersurfaces.
\bs\bs
\vbox{\centerline{\bf 1. Toric varieties}
\bsn
In this section we review some basics of toric varieties, mixed Hodge modules, intersection complexes, and prove Theorem~1.}
\msn
{\bf 1.1.~Toric varieties.} Let $\TT$ be an $n$-dimensional algebraic torus over $\C$. Let $M,N$ be the groups of characters and one-parameter subgroups of $\TT$ respectively, see \cite{Da}, \cite{KKMS}, \cite{Od}, etc. These are free abelian groups of rank $n$, and the composition induces a perfect pairing
$$M\times N\ni(u,v)\mapsto\langle u,v\rangle\in\Z.$$
\sk
Set $M_{\R}:=M\otimes_{\Z}\R$, and similarly for $N_{\R}$.
For a convex rational polyhedral cone $\eta$ in $N_{\R}$, we define the associated affine toric variety by
$$\aligned U_{\eta}&:={\rm Spec}\,\C[M\cap\eta^{\vee}]\q\q\q\q\h{with}\\ \eta^{\vee}&:=\bl\{u\in M_{\R}\,\big|\,\langle u,v\rangle\ges 0\,\,(\forall\,v\in\eta)\br\}.\endaligned
\leqno(1.1.1)$$
\sk
For a finite rational partial polyhedral decomposition (or a fan) $\Si$, we get the associated toric variety $\X_{\Si}$ by gluing the affine toric varieties $U_{\eta}$ ($\eta\in\Si$), see {\it loc.\,cit.}.
\sk
We have the canonical stratification by torus action orbits
$$\X_{\Si}=\msqcup_{\eta\in\Si\1\sqcup\1\{0\}}\,X_{\eta},
\leqno(1.1.2)$$
where
$$\aligned X_{\eta}&:={\rm Spec}\,\C[M\cap\eta^{\perp}]\q\q\q\q\h{with}\\ \eta^{\perp}&:=\bl\{u\in M_{\R}\,\big|\,\langle u,v\rangle= 0\,\,(\forall\,v\in\eta)\br\}.\endaligned
\leqno(1.1.3)$$
\msn
{\bf Remark~1.1.} The character of $\TT$ corresponding to $m\in M$ is denoted by $x^m$ in this paper. If we choose free generators $e_1,\dots,e_n$ of $M$, then $m$ and $x^m$ are respectively identified with $(m_1,\dots,m_n)\in\Z^n$ and $x_1^{m_1}\cdots x_n^{m_n}$, where $x_i$ is the character corresponding to $e_i$.
\msn
{\bf 1.2.~Dual fans.} Let $\De\subset M_{\R}$ be a polyhedron defined by a finite number of linear inequalities with rational coefficients. For each face $\si<\De$, we have the dual cone defined by
$$\aligned&\si^{\vee}:=\bl\{v\in\N_{\R}\,\big|\,\De_{v}\supset\si\br\}\q\q\h{with}\\&\De_{v}:=\bl\{u\in\De\,\big|\,\langle u,v\rangle=\min\langle\De,v\rangle\br\}.\endaligned
\leqno(1.2.1)$$
This gives the {\it dual fan\1} $\Si$ to $\De$, see for instance \cite{Va}. (It is also called the normal fan. However, min is usually replaced by max in the above definition, where the direction of the cone becomes opposite.) We then get the associated toric variety $\X_{\De}:=\X_{\Si}$.
We have the stratification by torus action orbits
$$\X_{\De}=\msqcup_{\si\les\De}\,X_{\si},
\leqno(1.2.2)$$
with $X_{\si}:=X_{\eta}$ for $\eta=\si^{\vee}$. Note that $\De^{\vee}=\{0\}$, and $X_{\De}=X_{\{0\}}=\TT$.
\msn
{\bf Remark~1.2.} If $\De$ is compact, there is a natural isomorphism
$$\X_{\De}={\rm Proj}\,\C[\Z^{n+1}{\cap}\1C_{\De}],
\leqno(1.2.3)$$
where $C_{\De}$ is as in a remark after Corollary~1. This seems to be used implicitly in \cite{BB}.
\msn
{\bf 1.3.~Quasi-smooth subdivisions.} We say that a fan $\Si$ is {\it quasi-smooth\1} if any cone $\eta\in\Si$ is simplicial, that is, spanned by $d_{\eta}$-lines with $d_{\eta}=\dim\eta$. In this case, it is known that the associated toric variety $\X_{\Si}$ is quasi-smooth, that is, a $V$-manifold. For any fan $\Si$, there is a ``standard" quasi-smooth subdivision $\Si'$ (which is combinatorially unique), see \cite[8.1]{Da}. The toric blow-up in the proof of Theorem~1 (see (1.7) below) gives the simplest examples when $n=3$ as in Example~(1.7) below.
\sk
If $\Si$ is the dual fan of a compact convex rational polyhedron $\De\subset M_{\R}$ (where rational means that the vertices are in $M_{\Q}$), this can be obtained as the dual fan of a {\it prime cutting\1} of $\De$ explained below. Here a compact polyhedron is called {\it prime\1} if its dual fan is simplicial.
\sk
Choose linear functions $\ell'_{\si}$ on $M_{\R}$ defined by $v_{\si}\in N_{\Q}$ (that is, $\ell_{\si}(u)=\langle u,v_{\si}\rangle$) together with rational numbers $a_{\si}$ for $\si<\De$ such that
$$\ell'_{\si}(\De)\subset\R_{\ges a_{\si}},\q\ell_{\si}^{\prime\,-1}(a_{\si})\cap\De=\si.$$
Set
$$\De':=\De\cap\mcap_{\si<\De}\,\ell_{\si}^{\prime\1-1}(\R_{\ges a_{\si}+\ep_{\si}}).
\leqno(1.3.1)$$
Here $\ep_{\si}\in\Q_{>0}$ with $\ep_{\tau}/\ep_{\si}$ ($\tau>\si$) {\it sufficiently small\1} (depending on the $\ell'_{\si}$) if the dual cone $\eta^{\vee}$ is {\it not\1} simplicial, and $\ep_{\si}=0$ otherwise (for instance, if $d_{\si}<3$). Note that the intersection over $\si<\De$ may be replaced by the one over $\si<\De$ with $\si^{\vee}$ {\it non-simplicial}. The dual fan of $\De'$ is a quasi-smooth subdivision explained above, where the $\ell'_{\si}$ correspond to the added 1-dimensional cones of $\Si'$.
(Note that a face of a simplicial cone is also simplicial.)
\sk
The above construction gives a prime polyhedron majorizing $\De$ in \cite{DK}. Indeed, one can proceed by induction on $n$ using a general hyperplane cut, since the hyperplanes $\ell^{\prime-1}_{\si}(a_{\si}{+}\ep_{\si})$ with $d_{\si}=0$ could be viewed as {\it general hyperplanes\1} if $\ep_{\tau}\ll\ep_{\si}\ll 1$ for $\tau>\si$. (Here a general hyperplane means that it intersects every positive-dimensional face transversally if the intersection is non-empty, and it does not intersect any vertices.)
\sk
When $\De$ is the Newton polyhedron of a non-degenerate Laurent polynomial, the $\ep_{\si}$ are considered to be {\it infinitesimally small\1} (or the ``limit" for $\ep_{\si}\to 0$ is considered). This gives a natural map $\pi$ from the set of faces of $\De'$ to that of $\De$.
\msn
{\bf Remark~1.3.} We can prove Theorem~1 except for the last assertion on the recursive formula by using the above prime cutting.
\msn
{\bf 1.4.~Some basics of mixed Hodge modules} (\cite{mhp}, \cite{mhm}, \cite{ypg}, etc.) For any complex algebraic variety $X$, there is an abelian category $\MHM(X)$ consisting of mixed Hodge modules on $X$. If $X$ is a point, there is a natural equivalence of categories
$$\MHM(pt)=\MHS,
\leqno(1.4.1)$$
where the right-hand side is the category of mixed $\Q$-Hodge structures (which are graded-polarizable, that is, the $\Gr^W_k$ are polarizable), see \cite[3.3.3]{ypg}.
\sk
Any mixed Hodge module $\M\in\MHM(X)$ has a canonical weight filtration $W$ such that the functor associating $\Gr^W_k\M$ to $\M$ is an exact functor, and every morphism of mixed Hodge modules is strictly compatible with $W$. We say that $\M$ has weights $\ges w$ (resp. $\les w$) if $\Gr^W_k\M=0$ for $k<w$ (resp. $k>w$). In particular, $\M$ is called pure of weight $w$ if $\Gr^W_k\M=0$ for $k\ne w$.
\sk
A pure Hodge module is {\it semi-simple,} and is a direct sum of simple pure Hodge modules with strict supports (using a strict support decomposition and a polarization). Here we say that $\M$ has {\it strict support\1} $Z$ if ${\rm supp}\,\M=Z$ and $\M$ has no non-trivial sub nor quotient objects in $\MHM(X)$ supported on a proper subset of $Z$. Moreover, a pure Hodge module with strict support $Z$ is generically a polarizable variation of Hodge structure on a smooth Zariski-open subset of $Z$, and conversely the latter can be extended uniquely to a pure Hodge module with strict support $Z$. Here we have some shifts of underlying $\Q$-complexes, see Remark~(1.4)\,(iii) below.
\sk
There is a forgetful functor
$${\rm rat}:D^b\MHM(X)\to D^b_c(X^{\rm an},\Q),$$
using the realization functor in \cite{BBD}, where ${\rm rat}(\M^{\ssb})$ is called the underlying $\Q$-complex of $\M^{\ssb}$.
The standard cohomological functor $H^j$ on the source corresponds to the cohomological functor ${}^p\Hc^j$ (see {\it loc.\,cit.}) on the target: there are canonical isomorphisms in $D^b_c(X^{\rm an},\Q)$
$${\rm rat}(H^j\M^{\ssb})={}^p\Hc^j{\rm rat}(\M^{\ssb})\q(j\in\Z).
\leqno(1.4.2)$$
The restriction of rat to $\MHM(X)\subset D^b\MHM(X)$ is faithful (that is, induces the injection of group of morphisms). These imply that $\M^{\ssb}\in D^b\MHM(X)$ is acyclic if ${\rm rat}(\M^{\ssb})$ is. Taking a mapping cone, we see that a morphism $u$ of $D^b\MHM(X)$ is an isomorphism is ${\rm rat}(u)$ is.
\sk
For any morphism of complex algebraic varieties $f:X\to Y$, there are canonically defined functors of bounded derived categories
$$\aligned f_*,f_!:D^b\MHM(X)&\to D^b\MHM(Y),\\ f^*,f^!:D^b\MHM(Y)&\to D^b\MHM(X),\endaligned$$
which are compatible with the corresponding functors of the underlying $\Q$-complexes via the functor rat, that is, there are canonical isomorphisms
$${\rm rat}(f_*\M^{\ssb})=\RR f_*\bl({\rm rat}(\M^{\ssb})\br),\,\,\h{etc}.
\leqno(1.4.3)$$
\sk
For complex algebraic varieties $X$, let $D^b\MHM(X)_{\les w}$, $D^b\MHM(X)_{\ges w}$, $D^b\MHM(X)_{(w)}$ be the full subcategories of $D^b\MHM(X)$ consisting of $\M^{\ssb}$ such that $H^j\M^{\ssb}$ has weights $\les w\pl j$, $\ges w\pl j$, and $w\pl j$ respectively. Then, for each integer $w$, we have
$$\aligned&\h{$D^b\MHM(X)_{\les w}$ are stable by $f_!,f^*$},\\ &\h{$D^b\MHM(X)_{\ges w}$ are stable by $f_*,f^*$},\\ &\h{$D^b\MHM(X)_{(w)}$ are stable by $f_!\eq f_*$ for $f$ proper.}\endaligned
\leqno(1.4.4)$$
see \cite[4.5.2]{mhm}. (These are analogous to \cite{BBD}.)
From the {\it semi-simplicity\1} of pure Hodge modules, we can deduce a non-canonical isomorphisms for $\M^{\ssb}\in D^b\MHM(X)_{(w)}$
$$\M^{\ssb}\cong\mopl_{j\in\Z}\,(H^j\M^{\ssb})[-j].
\leqno(1.4.5)$$
Combined with (1.4.4), this gives the {\it decomposition theorem\1} for pure Hodge modules. (An analogous assertion for $\ell$-adic sheaves seems quite non-trivial, since this seems to be related to the semi-simplicity of Frobenius action, although the decomposition theorem holds at least after forgetting Frobenius action, see \cite{BBD}. Note also that there does not seem to be a {\it canonical\1} increasing filtration on $\M^{\ssb}\in D^b\MHM(X)$ whose $k\1$th graded quotient is a pure complex of weight $k$ ($k\,{\in}\,\Z$). Indeed, these would depend heavily on the choice of representatives of $\M^{\ssb}$ in $C^b\MHM(X)$ if it is defined in a naive way.)
\sk
There are adjunction relations for $\M^{\ssb}\in D^b\MHM(X)$, $\Nc^{\ssb}\in D^b\MHM(Y)$ by the definition of pull-backs\,:
$$\aligned{\rm Hom}(f^*\Nc^{\ssb},\M^{\ssb})&={\rm Hom}(\Nc^{\ssb},f_*\M^{\ssb}),\\{\rm Hom}(\M^{\ssb},f^!\Nc^{\ssb})&={\rm Hom}(f_!\M^{\ssb},\Nc^{\ssb}).\endaligned$$
For a closed immersion $i:X\into Y$ and the complementary open immersion $j:Y\setminus X\into Y$, the adjunction relations induces some morphisms in the distinguished triangles
$$\aligned&j_!j^*\M^{\ssb}\to\M^{\ssb}\to i_*i^*\M^{\ssb}\buildrel{+1}\over\to,\\&i_*i^!\M^{\ssb}\to\M^{\ssb}\to j_*j^*\M^{\ssb}\buildrel{+1}\over\to,\endaligned
\leqno(1.4.6)$$
with $i_!=i_*$, $j^!=j^*$, see \cite[4.4.1]{mhm}. These triangles are dual of each other, and
$$\DD j_*=j_!\DD,\q \DD i^!=i^*\DD,
\leqno(1.4.7)$$
with $\DD$ the dual functor. This is a contravariant exact functor of $\MHM(X)$, and naturally induces a functor of $D^b\MHM(X)$ such that $\DD^2={\rm id}$.
\sk
If $X$ is an {\it affine\1} variety, then we have the {\it weak Lefschetz\1} type (or Artin's) theorem
$$\begin{array}{cl}H^j(X,\M):=H^j(a_X)_*\M=0&(j>0),\\ H_c^j(X,\M):=H^j(a_X)_!\M=0&(j<0),\raise5mm\h{}\end{array}
\leqno(1.4.8)$$
where $\M\in\MHM(X)$ with $a_X:X\to pt$ the structure morphism.
(This follows for instance from \cite[2.1.18]{mhp} using duality.)
\msn
{\bf Remark~1.4}\,(i). We can use (1.4.6) to show the equality
$$\chi_c(Y)=\chi_c(Y\setminus X)+\chi_c(X),
\leqno(1.4.9)$$
in the Grothendieck ring of mixed Hodge structures. Here
$$\chi(X):=\msum_j\,(-1)^j[H^j(X)],\q\chi_c(X):=\msum_j\,(-1)^j[H^j_c(X)],$$
with
$$H^j(X):=H^j(a_X)_*\Q_{h,X},\q H^j_c(X):=H^j(a_X)_!\Q_{h,X},$$
and $\Q_{h,X}:=a_X^*\Q_h$ with $\Q_h$ the pure Hodge structure of type $(0,0)$.
\msn
{\bf Remark~1.4}\,(ii). For a closed embedding $i:X\into Y$, the direct image $i_*$ is often omitted to simplify the notation. Indeed, there is an equivalence of categories
$$i_*:D^b\MHM(X)\simto D^b_X\MHM(Y),
\leqno(1.4.10)$$
where $D^b_X\MHM(Y)\subset D^b\MHM(Y)$ is the full subcategory of $\M^{\ssb}$ satisfying the condition that ${\rm supp}\,H^j\M^{\ssb}\subset X$, see \cite[4.2.10]{mhm}.
\msn
{\bf Remark~1.4}\,(iii). A polarizable variation of Hodge structure of weight $w$ on a smooth variety $X$ naturally defines a pure Hodge module $\M$ of weight $w\pl\dim X$ on $X$ such that ${\rm rat}(\M)$ is a local system {\it shifted by\1} $\dim X$ on $X$. Conversely, a pure Hodge module of weight $w\pl\dim X$ on a smooth variety $X$ is associated with a polarizable variation of Hodge structure of weight $w$ as above if ${\rm rat}(\M)$ is a local system on $X$ shifted by $\dim X$.
\msn
{\bf 1.5.~Intersection complexes.} Let $X$ be an irreducible complex algebraic variety. Let $X^{\rm an}$ be its associated analytic space. Set $n:=\dim X$. The {\it intersection complex\1} of $X^{\rm an}$ can be defined in $D^b_c(X^{\rm an},\Q)$ by
$${\rm IC}_{X}\Q=\tau_{<0}\RR(j_0)_*\cdots\tau_{<1-n}\RR(j_{n-1})_*(\Q_{U^{\rm an}_n}[n]),
\leqno(1.5.1)$$
where the $j_i:U_{i+1}\into U_i$ denote the inclusions with $\{U_i\}$ a decreasing sequence of Zariski-open subsets of $X$ such that the $S_i:=U_i\setminus U_{i+1}$ are smooth of pure dimension $i$, and form a {\it Whitney stratification,} see \cite{BBD}. This is independent of the choice of $\{U_i\}$ as long as Whitney's condition~(b) is satisfied for $\{S_i\}$.
\sk
The last condition is satisfied for the case of the natural stratification in (1.2.2) with $X=\X_{\De}$, since the singularities of $\X_{\De}$ are {\it constant\1} along each stratum $X_{\si}$ by definition. This also implies that
$$\Hc^j{\rm IC}_{\X_{\De}}\Q|_{X_{\si}}\,\,\,\h{is a constant sheaf for any}\,\,j,\,\si.
\leqno(1.5.2)$$
\sk
The {\it intersection cohomology\1} of an irreducible complex algebraic variety $X$ is defined by
$${\rm IH}^j(X):=H^j(X^{\rm an},{\rm IC}^s_X\Q)\q\h{with}\q{\rm IC}^s_X\Q:={\rm IC}_X\Q[-n].
\leqno(1.5.3)$$
We have similarly the intersection cohomology with compact supports ${\rm IH}^j_c(X)$ by replacing $H^j$ with $H^j_c$. These have canonical mixed Hodge structures by using ${\rm IC}_Z\Q_h$, which is a pure Hodge module of weight $n$ such that its underlying $\Q$-complex is ${\rm IC}_Z\Q$, where $n=\dim Z$. It is also characterized by the condition that it has strict support $Z$ and its restriction to a dense Zariski-open smooth subset $U\subset Z$ is $\Q_{h,U}[\dim U]$.
\sk
We have the duality isomorphism
$$\DD\1{\rm IH}^j(X)={\rm IH}_c^{2n-j}(X)(n),
\leqno(1.5.4)$$
where $\DD\1{\rm IH}^j(X):={\rm Hom}_{\Q}\bl({\rm IH}^j(X),\Q\br)$ which is the dual mixed Hodge structure.
\sk
If $Y$ is an irreducible projective variety with $\ell$ the first Chern class of an ample line bundle, we have by \cite{mhp} the {\it purity} :
$${\rm IH}^j(Y)\,\,\,\h{is pure of weight}\,\,\,j\,\,\,(\forall\,j\in[0,2n]),
\leqno(1.5.5)$$
together with the {\it hard Lefschetz\1} property\,:
$$\ell^j:{\rm IH}^{n-j}(Y)\simto{\rm IH}^{n+j}(Y)(j)\q(\forall j>0).
\leqno(1.5.6)$$
\msn
{\bf 1.6.~Generalized Thom-Gysin sequence.} Let $X'$ be an algebraic $\C^*$-bundle over an irreducible algebraic variety $Y$. We have a {\it generalized Thom-Gysin sequence\1}
$$\to{\rm IH}^{j-1}(X')\to{\rm IH}^{j-2}(Y)(-1)\buildrel{\xi}\over\to{\rm IH}^j(Y)\to{\rm IH}^j(X')\to
\leqno(1.6.1)$$
in the category of mixed Hodge structures, where the morphism $\xi$ is given by the action of the first Chern class of the line bundle associated with the $\C^*$-bundle up to a sign.
\sk
This can be proved by considering a line bundle $\pi:Z\to Y$ such that the complement of the zero section is $X'$, and using the long exact sequence of mixed Hodge structures associated with the direct image by $\pi_*$ of the following distinguished triangle in $D^b\MHM(Z)$\,:
$$(i_Y)_*i_Y^!{\rm IC}^s_Z\Q_h\to{\rm IC}^s_Z\Q_h\to(j_{X'})_*{\rm IC}^s_{X'}\Q_h\buildrel{+1}\over\to.
\leqno(1.6.2)$$
Here $Y$ is identified with the zero section of Z, and $i_Y:Y\into Z$, $j_{X'}:X'\into Z$ are natural inclusions. We have by definition
$${\rm IC}^s_Z\Q_h:={\rm IC}_Z\Q_h[-n].$$
There are canonical isomorphisms in $D^b\MHM(Y)$\,:
$$i_Y^!{\rm IC}^s_Z\Q_h={\rm IC}^s_Y\Q_h(-1)[-2],\q\pi_*{\rm IC}^s_Z\Q_h={\rm IC}^s_Y\Q_h.
\leqno(1.6.3)$$
The calculation of the morphism $\xi$ is reduced to the constant sheaf case, see \cite[1.3]{RSW}.
\sk
In the case $X'$ is the complement of the zero section of an ample line bundle on a projective variety $Y$, the generalized Thom-Gysin sequence (1.6.1) implies the isomorphisms of mixed Hodge structures
$${\rm IH}^j(X')=\begin{cases}{\rm IH}^j_{\rm prim}(Y)&(j<n),\\{\rm IH}^{j-1}_{\rm coprim}(Y)(-1)&(j\ges n),\end{cases}
\leqno(1.6.4)$$
using the {\it primitive decomposition\1} associated with the hard Lefschetz property (1.5.6)\,:
$$\mopl_{j=0}^{2n-2}{\rm IH}^j(Y)=\mopl_{i=0}^{n-1}\,\mopl_{k=0}^i\,\ell^k\1{\rm IH}^{n-1-i}_{\rm prim}(Y)(-k),
\leqno(1.6.5)$$
(Note that $\dim Y=n\mi 1$.) Here we have for $k\ges 0$
$$\aligned{\rm IH}^{n-1-k}_{\rm prim}(Y)&:={\rm Ker}\,\ell^{k+1}\subset{\rm IH}^{n-1-k}(Y),\\{\rm IH}^{n-1+k}_{\rm coprim}(Y)&:={\rm Ker}\,\ell\subset{\rm IH}^{n-1+k}(Y),\endaligned
\leqno(1.6.6)$$
and these are 0 if $k<0$. There are isomorphisms for $k>0$
$$\ell^k:{\rm IH}^{n-1-k}_{\rm prim}(Y)\simto{\rm IH}^{n-1+k}_{\rm coprim}(Y)(k).
\leqno(1.6.7).$$
\msn
{\bf Remark~1.6.} By (1.6.4), the intersection cohomology ${\rm IH}^j(X')$ is pure of weight $j$ if $j\1{<}\1n$, and it is pure of weight $j\pl1$ otherwise. Both assertions are needed to show that the truncation by degree can be done by the truncation by weight.
Indeed, there may be {\it cancelations\1} during the passage to the Grothendieck ring if the above conditions are not satisfied (or if the hard Lefschetz theorem does not hold).
\msn
{\bf 1.7.~Proof of Theorem~1.} We may assume that the toric variety $\X_{\De}$ is {\it affine,} that is, $\De=\eta^{\vee}$ for the unique $n$-dimensional cone $\eta$ in $\Si$. Choose a 1-dimensional rational cone $\eta_0$ inside $\eta=\De^{\vee}$. We have the subdivision $\Si'$ of $\Si$ by taking the convex hull of $\eta_0\cup\eta'$ for any $\eta'\in\Si$ with $\dim\eta'<n$. This can be obtained also as the dual fan of $\De'\subset\De$ defined by
$$\De':=\bl\{u\in\De\,\big|\,\langle u,v\rangle>c\1\br\},$$
where $v\in(N\cap\eta_0)\setminus\{0\}$ and $c\in\Q_{>0}$. Set
$$\De'':=\bl\{u\in\De'\,\big|\,\langle u,v\rangle=c\1\br\}.$$
This is the unique $(n{-}1)$-dimensional compact face of $\De'$.
\sk
Let $\X_{\De'}$ be the toric variety associated with $\De'$, that is, with $\Si'$. We have the projective morphism
$$\pi:\X_{\De'}\to\X_{\De},$$
inducing an isomorphism outside the vertex $0$ of the affine cone $\X_{\De}$. This is a toric blow-up of $\X_{\De}$. Set
$$A_{\De}:=\C[M\cap\De].$$
This is the affine ring of $\X_{\De}$. We have the isomorphism
$$\X_{\De''}={\rm Proj}\,A_{\De}.$$
Here the left-hand side is defined by taking the $(n{-}1)$-dimensional affine space containing $\De''$, and the grading of $A_{\De}$ is defined by the pairing with the primitive element $v\in(N\cap\eta_0)\setminus\{0\}$. The exceptional divisor $\pi^{-1}(0)$ is identified with these varieties. Moreover $\X_{\De'}$ is a line bundle over $\X_{\De''}={\rm Proj}\,A_{\De}$, and the restrictions of the fibers to $\X_{\De}\setminus\{0\}$ are identified with the orbits of the action of the one-parameter subgroup corresponding to $v$. (Indeed, $\C[M\cap \ell_v^{-1}(0)]$ is identified with functions on $\TT$ which are constant on each orbit of the $\C^*$-action, where $\ell_v$ is the linear function on $M_{\R}$ defined by the pairing with the above $v$.)
\sk
To simplify the notation, set
$$Y:=\X_{\De''}={\rm Proj}\,A_{\De},\q X:=\X_{\De},\q Z:=\X_{\De'},\q X':=X\setminus\{0\}=Z\setminus Y,$$
so that $X'$ is a $\C^*$-bundle over $Y$ as in (1.6). By the decomposition theorem (see (1.4.4--5)), there is a non-canonical isomorphism in $D^b{\rm MHM(X)}$\,:
$$\pi_*{\rm IC}_Z\Q_h\cong{\rm IC}_Z\Q_h\oplus\mopl_{j\in\Z}\,M_0^j[n{-}j],
\leqno(1.7.1)$$
with $M_0^j$ pure Hodge modules of weight $j$ on $\{0\}$, which are identified with pure Hodge structures of the same weight. (Note that the category of mixed Hodge modules on $pt$ is identified with that of mixed Hodge structures, see (1.4.1).) More precisely, $M_0^j$ is defined to be the direct factor of $H^{j-n}\pi_*{\rm IC}_Z\Q_h$ supported on $\{0\}$ (which is equal to $H^{j-n}\pi_*{\rm IC}_Z\Q_h$ if $j\ne n$). The direct image by the inclusion $i_0:\{0\}\into X$ is omitted to simplify the notation, see Remark~(1.4)\,(ii).
\sk
There are isomorphisms of pure Hodge structures
$$\ell^j:M_0^{n-j}\simto M_0^{n+j}(j)\q(j>0),
\leqno(1.7.2)$$
by the hard Lefschetz theorem for the direct image of pure Hodge modules by projective morphisms \cite[Theorem~1]{mhp}. Here $\ell$ is the first Chern class of the relatively ample line bundle of $\pi$, and its restriction to $Y$ coincides with that for an ample line bundle, since $Y$ is the exceptional divisor. (Note that the action of $\ell$ {\it in the derived category\1} is not compatible with the direct sum decomposition of the right-hand side of (1.7.1), and there may be elements in the higher extension groups between the direct factors of the right-hand side.)
\sk
Taking the pull-back by the inclusion $i_0:\{0\}\into X$, and applying the base change theorem (see \cite[(4.4.3)]{mhm}), we can get non-canonical isomorphisms of pure Hodge structures of weight $j$\,:
$${\rm IH}^j(Y)\cong H^ji_0^*{\rm IC}^s_X\oplus M_0^j\q(j\in\Z),
\leqno(1.7.3)$$
since $i_Y^*{\rm IC}^s_Z\Q_h={\rm IC}^s_Y\Q_h$. By the hard Lefschetz theorem (1.5.6), there are isomorphisms
$$\ell^j:{\rm IH}^{n-1-j}(Y)\simto{\rm IH}^{n-1+j}(Y)\q(j>0).
\leqno(1.7.4)$$
Using the {\it primitive decompositions\1} as in (1.6.5) together with the {\it Grothendieck ring\1} of mixed Hodge structures as well as the purity of the ${\rm IH}^j(Y)$ and also the semisimplicity of polarizable pure Hodge structures, we then conclude that $\mopl_j\,M_0^j$ is non-canonically isomorphic to the {\it non-primitive part\1} of $\mopl_j\,{\rm IH}^j(Y)$, and moreover there are non-canonical isomorphisms of pure Hodge structures
$$H^ji_0^*{\rm IC}^s_X\Q_h\cong{\rm IH}^j_{\rm prim}(Y)\q(\forall\,j\in\Z),
\leqno(1.7.5)$$
since $[H^ji_0^*{\rm IC}^s_X\Q_h]=[{\rm IH}^j_{\rm prim}(Y)]$.
Note that ${\rm IH}^j_{\rm prim}(Y)=0$ if $j\ges n$.
\sk
By the purity (1.5.5) together with (1.7.5), we can calculate the {\it Poincar\'e polynomials\1} of ${\rm IH}^{\ssb}(Y)$, $H^{\ssb}i_0^*{\rm IC}^s_X\Q_h$ using their classes in the Grothendieck ring of mixed Hodge structures $K_0(\MHS)$, more precisely, using
$$\aligned{[{\rm IH}^{\ssb}(Y)]}&:=\msum_j\,(-1)^j[{\rm IH}^j(Y)],\\ [i_0^*{\rm IC}^s_X\Q_h]&:=\msum_j\,(-1)^j[H^ji_0^*{\rm IC}^s_X\Q_h],\endaligned$$
since {\it the degree is uniquely determined by the weight}.
\sk
Let $K_0(\MHS)^{\rm Tate}$ be the subring of $K_0(\MHS)$ which is generated by non-negative powers of $t:=[\Q(-1)]$. This is identified with the polynomial ring $\Z[t]$. We first show that
$$[{\rm IH}^{\ssb}(Y)],\,[i_0^*{\rm IC}^s_X\Q_h]\,\in\,K_0(\MHS)^{\rm Tate}.
\leqno(1.7.6)$$
This implies that ${\rm IH}^{\ssb}(Y)$ and $H^{\ssb}i_0^*{\rm IC}^s_X\Q_h$ have only even weights, hence the odd degree part vanishes by the purity.
\sk
We have the stratification as in (1.2.2)
$$Y=\msqcup_{\tau\les\De''}\,X_{\tau}.$$
There is a one-to-one correspondence between the faces $\tau\les\De''$ and $\si\les\De$ except for the origin $0\les\De$ by considering the cone $C(\tau)$ or taking the intersection $\si\cap\De''$ inversely.
Since $X'$ is a $\C^*$-bundle over $Y$, we have the equalities
$$m_{\si,j}=m_{\tau,j}\,\,\,(j\in\Z)\q\h{if}\q\si=C(\tau),
\leqno(1.7.7)$$
in the notation of (2). (This implies that the formulation in \cite{St} is equivalent to the one in \cite{DL}, \cite{Fi}.)
Moreover we have for $\tau\les\De''$
$$[H_c^{\ssb}(X_{\tau})]=(t\mi1)^{\dim\tau}\,\in\,K_0(\MHS)^{\rm Tate},
\leqno(1.7.8)$$
since $X_{\tau}\cong(\C^*)^{\dim\tau}$. Then (1.7.6) follows from (1.7.5) inductively taking the intersection with a rationally defined affine subspace intersecting $\si$ transversely for each $0<\si\les\De$.
This proves Theorem~1 except for the last assertion by using (1.5.2), since a mixed Hodge module on a smooth variety such that its underlying $\Q$-complex is a local system is a variation of mixed Hodge structure (up to a shift of complex).
\sk
For the proof of the last assertion of Theorem~1, set for $\si\les\De$
$$m_{\si}(t):=\msum_{k\in\N}\,m_{\si,2k}\1t^k\,\in\,\Z[t].$$
The above argument implies the recursive formula
$$m_0(t)=\tau_{<n/2}\bl((1\mi t)\1\msum_{0<\si\les\De}\,(t\mi 1)^{\dim\si-1}m_{\si}(t)\br),
\leqno(1.7.9)$$
using (1.7.5). Here $\tau_{<\al}\,h(t):=\msum_{k<\al}\,c_k\1t^k$ if $h(t)=\msum_{k\in\N}\,c_k\1t^k$ $(c_k\in\Z,\,\al\in\Q$). Indeed, we have
$$[{\rm IH}^{\ssb}(Y)]=\msum_{0<\si\les\De}\,(t\mi 1)^{\dim\si-1}m_{\si}(t),
\leqno(1.7.10)$$
using (1.7.7--8). The truncation $\tau_{<n/2}$ after the multiplication by $(1\mi t)$ gives the passage to the primitive part. This finishes the proof of Theorem~1.
\msn
{\bf Remark~1.7.} The isomorphisms in (1.7.5) are closely related to (1.6).
Indeed, (1.5.1) implies the canonical isomorphisms
$$H^ji_0^*{\rm IC}^s_X\Q_h\simto\begin{cases}H^ji_0^*(j_0)_*{\rm IC}^s_{X'}\Q_h&(j<n),\\ \,0&(j\ge n).\end{cases}
\leqno(1.7.11)$$
with $j_0:X'\into X$ the natural inclusion. Note that we have the canonical morphisms induced by ${\rm id}\to(j_0)_*j_0^*$, and the isomorphisms hold for the underlying $\Q$-complexes by (1.5.1). Similarly there are canonical isomorphisms
$$H^j(X',{\rm IC}^s_{X'}\Q_h)\simto H^ji_0^*(j_0)_*{\rm IC}^s_{X'}\Q_h,
\leqno(1.7.12)$$
see (1.4.8) for the left-hand side. The canonical morphisms are induced by ${\rm id}\to(i_0)_*i_0^*$. They are isomorphisms for the underlying $\Q$-complexes using the $\C^*$-action. Note that the canonical mixed Hodge structure on ${\rm IH}^j(X')$ is defined by the left-hand side of (1.7.12). Combined with (1.6.4), we thus get the isomorphisms of mixed Hodge structures (1.7.5).
\sk
Note also that there are equalities
$$\aligned\h{$[{\rm IH}^{\ssb}(X')]$}&=(1\mi t)\1\msum_{0<\si\les\De}\,(t\mi 1)^{\dim\si-1}m_{\si}(t),\\
[{\rm IH}^{\ssb}_c(X')]&=\msum_{0<\si\les\De}\,(t\mi 1)^{\dim\si}m_{\si}(t),\endaligned
\leqno(1.7.13)$$
respectively by (1.7.10), (1.6.1) (or a spectral sequence) and by (1.7.8) with $\tau$ replaced by $\si$. In particular, we get the equality
$$[{\rm IH}^{\ssb}(X')]=-[{\rm IH}_c^{\ssb}(X')].$$
This follows also from the property that $X'$ is a $\C^*$-bundle over a proper variety $Y$ (using some spectral sequences). However, the truncation of ${\rm IH}^{\ssb}(X')$ by degree $<n/2$ cannot be replaced with that of $[{\rm IH}^{\ssb}(X')]$ by weight $<n$ {\it without\1} using (1.6.4) which is a consequence of the purity, the hard Lefschetz theorem, and the generalized Thom-Gysin sequence.
\msn
{\bf Example~1.7.} Assume $n=3$, and $\X$ is an affine toric variety as in the proof of Theorem~1. Let $k_0$ is the number of 2-dimensional faces of $\De\subset\R^3$ passing through the vertex $0\in\De$ (which can be identified with a point of $\X$). Then
$$\aligned&m_0(t)=\tau_{<3/2}(1\mi t)\bl((t\mi 1)^2+k_0(t\mi 1)+k_0\br)=1+(k_0\mi 3)t,\\&\h{that is,}\q\q\q\dim_{\Q}\Hc^{-1}({\rm IC}_{\X}\Q)_0=k_0\mi 3.\endaligned$$
We have $\chi(Y)=k_0$, and
$$\dim H^k(Y)=\begin{cases}k_0\mi 2&(j=2),\\ 1&(j=0,4),\\0&\h{(otherwise).}\end{cases}$$
Hence
$$\RR\pi_*{\rm IC}_Z\Q\cong{\rm IC}_X\Q\oplus\Q_0(-2)[-1]\oplus\Q_0(-1)[1].$$
\bs\bs
\vbox{\centerline{\bf 2. Toric hypersurfaces}
\bsn
In this section, we explain some applications of Theorem~1 to the study of the cohomology of non-degenerate toric hypersurfaces.}
\msn
{\bf 2.1.~Non-degenerate toric hypersurfaces.} Let $f$ be a section of $\OO_{\TT}$. This is identified with a Laurent polynomial by choosing free generators of $M$ as above. We have
$$f=\msum_{m\in M}\,f_m\1x^m\q(f_m\in\C),$$
where $f_m=0$ except for a finite number of $m$. The {\it support\1} of $f$ (denoted by ${\rm Supp}\,f$) is the union of $m\in M$ with $f_m\ne 0$. The {\it Newton polyhedron\1} $\Ga(f)$ of $f$ is the convex hull of the support of $f$ in $M_{\R}$.
\sk
For a face $\si\les\Ga(f)$, set
$$f_{\si}:=\msum_{m\in{\si}}\,f_m\1x^m.$$
We say that $f$ is {\it Newton non-degenerate\1} if $f_{\si}^{-1}(0)\subset\TT$ is {\it smooth\1} (in particular, reduced) for any $\si\les\Ga(f)$ including the case $\si=\Ga(f)$.
\ms
The above condition on each $f_{\si}$ is equivalent to the one in \cite{Ko} (using the partial derivatives of $f_{\si}$) if the minimal affine space containing $\si$ does not pass through the origin. Note that the condition is not imposed for $\si=\Ga(f)$ in {\it loc.\,cit.}
\ms
Let $\Zof$ be the closure of $\Zf:=f^{-1}(0)\subset\TT$ in $\X:=\X_{\De}$ with $\De=\Ga(f)$. Here $f$ is a Newton non-degenerate Laurent polynomial, and we assume $n\ges 2$. We call $\Zof$ (and also $\Zf$) a {\it non-degenerate toric hypersurface\1} with respect to $\De$. We have the stratification induced by (1.2.2)
$$\Zof=\msqcup_{\si\les\De}\,\Zfs\q\h{with}\q\Zfs:=\Zof\cap X_{\si}
\leqno(2.1.1)$$
\sk
The hypersurface $\Zof\subset\X$ is {\it stratified-smooth}. This means that it is locally defined by a function $h$ on a smooth ambient variety containing $\X$ such that the restriction of $\ddd h$ to each stratum $X_{\si}$ with $d_{\si}>0$ does not vanish on a neighborhood of $\Zfs$. This shows that the singularities of $\Zof$ is essentially the same as those of $\X$ locally. (Note that these are {\it constant\1} along each stratum.) It implies also that the union of $D:=\X\setminus\TT$ and $\Zof$ is a {\it divisor with normal crossings\1} if so is $D$ with $\X$ smooth.
\ms
By the above argument we get the canonical isomorphism in the notation of (1.5--6)
$$i_{\Zof}^*{\rm IC}^s_{\X}\Q_h={\rm IC}^s_{\Zof}\Q_h,
\leqno(2.1.2)$$
with $i_{\Zof}:\Zof\into\X$ the natural inclusion. Indeed, $H^{-1}i_{\Zof}^*{\rm IC}_{\X}\Q_h$ is a pure Hodge module of weight $n{-}1$, since this is a non-characteristic restriction with codimension 1. Here non-characteristic means that $\varphi_h$ vanishes with $h$ a local defining function of $\Zof\subset\X$.
\msn
{\bf 2.2.~Proof of Corollary~2.} Consider the distinguished triangle
$$\M\to\Q_{\Zof,h}\to{\rm IC}^s_{\Zof}\Q_h\buildrel{+1}\over\to,
\leqno(2.2.1)$$
in the notation of (1.4--6). The second morphism of (2.2.1) follows from \cite[(4.5.11)]{mhm}, and $\M$ is its mapping cone up to a shift of complex.
\sk
For $x\in\Zfs$, let $i_x:\{x\}\into\Zof$ be the inclusion. Theorem~1 together with (2.1.2) implies that
$$H^ji_x^*\M=\begin{cases}\Q_h^{m_{\si,j-1}}(-k)&\h{if}\,\,\,j=2k{+}1,\,k\in\Z_{>0},\\ \,0&\h{otherwise}.\end{cases}
\leqno(2.2.2)$$
Corollary~2 then follows from (2.2.1--2) together with (1.5.2) and (2.1.2). Indeed, $k$ in (2.2.2) is strictly positive, and the Hodge level of the following pure Hodge structure is at most $w\mi 2k$\,:
$$\Gr^W_w\bl(H^{i-2k-1}_c(\Zfs)(-k)\br)=\bl(\Gr^W_{w-2k}H^{i-2k-1}_c(\Zfs)\br)(-k).$$
Indeed, $\Gr^W_wH^i_c(\Zfs)$ has Hodge level $\les w$ independently of $i$ as is well-known. (Recall that the Hodge level of a pure Hodge structure is the difference between the maximum and the minimum of the integers $p$ such that $\Gr_F^p$ does not vanish.) This finishes the proof of Corollary~2.
\msn
{\bf 2.3.~Construction of \v Cech type complexes.} Set $V:=M_{\R}$. Let $\De\subset V$ be a polyhedron of dimension $n$ defined by a finite number of inequalities, and $\De^{\circ}$ be the interior of $\De$. Denote by $j_S:S\into V$ the inclusion of a locally closed subset $S\subset V$ in general. We have the following.
\msn
{\bf Proposition~2.3.} {\it There is a quasi-isomorphism of $\Z$-complexes on $\De$\,:
$$(j_{\De^{\circ}})_!\Z_{\De^{\circ}}\simto C_{\De}^{\1\ssb}\q\h{with}\q C_{\De}^{\1p}:=\mopl_{c_{\si}=p}\,\Z_{\si}.
\leqno(2.3.1)$$
Here the direct sum is taken over the faces $\si\les\De$ with $c_{\si}=p$, and $\Z_{\De^{\circ}}$, $\Z_{\si}$ are constant sheaves on $\De^{\circ}$, $\si$ respectively with stalks $\Z$. The direct image by closed embeddings are omitted to simplify the notation.}
\msn
{\it Proof.} Consider a finite increasing filtration $G$ on the constant sheaf $\Z_{\De}$ such that
$$G_p\1\Z_{\De}=(j_{\De_{(p)}})_!\1\Z_{\De_{(p)}},\q\Gr_p^G\1\Z_{\De}=\mopl_{c_{\si}=p}\,(j_{\si^{\circ}})_!\1\Z_{\si^{\circ}},$$
where $\De_{(p)}:=\msqcup_{c_{\si}\les p}\,\si^{\circ}$.
We have the isomorphism
$$\Z_{\De}=\RR(j_{\De^{\circ}})_*\1\Z_{\De^{\circ}},$$
and similarly for $\si$.
\sk
For an orientable {\it real\1} manifold $S$ of dimension $d_S$., the dual functor $\DD$ can be defined by
$$\DD\1K^{\ssb}:=\RR\Hc om_{\Z}(K^{\ssb},\Z_S[d_S])\q\bl(K^{\ssb}\in D^b_c(S,\Z)\br),$$
choosing an orientation of $S$. For locally closed real submanifold $S\subset V$, there is a functorial isomorphism \cite{Ve}
$$\DD\ssc(j_S)_!\cong\RR(j_S)_*\ssc\DD.$$
(This implies that $\DD$ commutes with the direct images by closed immersions.) Combined with the above isomorphism with $\De$ replaced by $\si$, it gives the isomorphisms
$$\DD\1(j_{\si^{\circ}})_!\1\Z_{\si^{\circ}}\cong\Z_{\si}[d_{\si}]\q(\si\les\De),
\leqno(2.3.2)$$
choosing an orientation of the real manifold $\si^{\circ}$. We have also the isomorphism
$$\DD\1\Z_{\De}=(j_{\De^{\circ}})_!\Z_{\De^{\circ}}[n].
\leqno(2.3.3)$$
This follows from the so-called Verdier {\it biduality\1} $\DD^2\1\cong{\rm id}$. (In this case, we can use also the distinguished triangles
$$(j_{V\setminus\De})_!\1\Z_{V\setminus\De}\to\Z_V\to\Z_{\De}\buildrel{+1}\over\to$$
$$(j_{\De^{\circ}})_!\1\Z_{\De^{\circ}}\to\Z_V\to\Z_{V\setminus\De^{\circ}}\buildrel{+1}\over\to$$
comparing the dual of the first with the second shifted by $n$.)
\sk
Applying $\DD$ to the filtered sheaf $(\Z_{\De},G)$, and shifting the complex by $-n$, we then get a finite decreasing filtration $G$ on
$$(\DD\1\Z_{\De})[-n]\cong(j_{\De^{\circ}})_!\Z_{\De^{\circ}},$$
satisfying
$$\Gr^p_G(j_{\De^{\circ}})_!\Z_{\De^{\circ}}\cong\mopl_{c_{\si}=p}\,\Z_{\si}[-c_{\si}]\q(p\in\Z).
\leqno(2.3.4)$$
We may now apply \cite[3.1.8--9]{BBD} (with respect to the {\it classical\1} $t$-structure). This gives the above complex $C_{\De}^{\1\ssb}$ together with the desired quasi-isomorphism, choosing an {\it orientation\1} of each $\si\les\De$. Note that the differential of $C_{\De}^{\1\ssb}$ is given by linear combinations of restriction morphisms with coefficients 1 or $-1$ which depend on the choice of orientations of the $\si\les\De$.
This finishes the proof of Proposition~(2.3).
\msn
{\bf Corollary~2.3.} {\it There are quasi-isomorphism of $\Z$-complexes on $\X$
$$\aligned(j_{X_{\De}})_!\Z_{X_{\De}}\simto K_{\X}^{\ssb}\q&\h{with}\q K_{\X}^p:=\mopl_{c_{\si}=p}\,\Z_{\Xo_{\si}},\\(j_{\Zf})_!\Z_{\Zf}\simto K_{\Zof}^{\ssb}\q&\h{with}\q K_{\Zof}^p:=\mopl_{c_{\si}=p}\,\Z_{\Zofs},\endaligned
\leqno(2.3.5)$$
where $j_{X_{\De}}:X_{\De}\,(=\TT)\into\X$, $j_{\Zf}:\Zf\into\Zof$ are natural inclusions.}
\msn
{\it Proof.} Since the combinatorial data of the faces $\si\les\De$ and the intersections of irreducible components of $\X\setminus\TT$ (that is, the $\Xo_{\si}$ for $\si\les\De$) are the same, the first assertion follows from Proposition~(2.3). (Indeed, the stalk of $K_{\X}^{\ssb}$ at $x\in X_{\si}$ can be identified with that of $C_{\De}^{\ssb}$ at $v\in\si^{\circ}$.)
This implies the second assertion by restricting to $\Zof\subset\X$. This finishes the proof of Corollary~(2.3).
\msn
{\bf Remark~2.3.} Proposition~(2.3) implies the following equalities for any $v\in\dd\1\De$\,:
$$\msum_{p=0}^n\,(-1)^pn_{\De,v,p}=0,
\leqno(2.3.6)$$
where
\vskip-2mm\nin
$$n_{\De,v,p}:=\#\{\1\si\les\De\mid\si\ni v,\,c_{\si}=p\1\}.$$
\msn
{\bf 2.4.~Calculation of certain Hodge numbers.} The formula (1) in the introduction (see also Remark~(2.4) below) can be made more precise as follows:
$$h^{k,p,q}_c(Z)=\begin{cases}\tbinom{n}{i}&\bl(k\eq 2n{-}2{-}i,\,p\eq q\eq n{-}1{-}i\,\,(i\in[0,n{-}2])\br),\\ \,\,\,0&\bl(k\,{<}\,n{-}1\,\,\,\h{or}\,\,\,k\eq n{-1}\,{<}\,p{+}q\,\,\,\h{or}\,\,\,k\,{>}\,n{-}1,p\,{\ne}\,q\br).\raise5mm\h{}\end{cases}
\leqno(2.4.1)$$
So there remains only the case $k\eq n{-}1\,{\ges}\,p{+}q$.
\sk
The following seems to be known to specialists:
\msn
{\bf Proposition~2.4.} {\it There are equalities
$$h^{n-1,n-1,0}_c(\Zf)=h^{n-1,n-1,0}(\Zof)=\bl|\De^{\circ}\cap M\br|,
\leqno(2.4.2)$$
where $\De^{\circ}$ is the interior of $\De$.}
\msn
{\it Proof.} The first equality follows from a long exact sequence of mixed Hodge structures, since $\dim\Zof\setminus\Zf<n{-}1$.
The last equality of (2.4.2) follows from the assertion that the reduced ideal $\I_D\subset\OO_{\X}$ of the boundary $D:=\X\setminus\TT$ can be identified with the dualizing sheaf $\om_{\X}$ of $\X$ using the freeness of the logarithmic differential forms $\Om_{\X}^{\1n}(\log D)$ (without using the calculation of interior Ehrhart polynomials). Indeed, there is a short exact sequence
$$0\to\om_{\X}\to\om_{\X}(\Zof)\to\om_{\Zof}\to 0,
\leqno(2.4.3)$$
since $\Zof$ is stratified-smooth. Here $H^q(\X,\om_{\X})=0$ ($q<n)$ by Grothendieck duality. Note that toric varieties are normal and Cohen-Macaulay, and moreover they have at worst {\it rational singularities\1} using the above freeness of logarithmic differential forms and studying the above ideal sheaf under a desingularization by a smooth subdivision of the fan. Since $\Zof$ is stratified smooth, this also implies that $\om_{\Zof}$ coincides with the direct image of the dualizing sheaf of a desingularization of $\Zof$, where the higher direct images vanish by Grauert-Riemenschneider vanishing theorem. We then get
$$h^{n-1,n-1,0}(\Zof)=\dim H^0(\Zof,\om_{\Zof}),
\leqno(2.4.4)$$
since the first equality of (2.4.2) holds with $\Zof$ replaced by any compactification of $\Zf$.
\sk
It is also known that
$$\dim H^0\bl(\X,\I_D\otimes_{\OO}\OO_{\X}(\Zof)\br)=\bl|\De^{\circ}\cap M\br|.
\leqno(2.4.5)$$
since the linear system associated with $\Zof$ corresponds to the line bundle associated with $\De$. So the last equality of (2.4.2) follows from the long exact sequence associated with (2.4.3). This finishes the proof of Proposition~(2.4).
\msn
{\bf Remark~2.4.} It is not necessarily trivial to show (1) in the introduction. For the assertion in the case $j>n{-}1$, we need the {\it weak Lefschetz\1} type (or Artin's) theorem (1.4.8) together with the ampleness of $\Zof\subset\X$ and also the canonical isomorphism of $\Q$-complexes
$$i_{\Zof}^!(j_{\TT})_!\1\Q_{\TT}=(j_{\Zf})_!\1\Q_{\Zf}(-1)[-2],
\leqno(2.4.6)$$
with $i_{\Zof}:\Zof\into\X$, $j_{\TT}:\TT\into\X$, $j_{\Zf}:\Zf\into\Zof$ natural inclusions. The latter isomorphism follows from the {\it stratified smoothness\1} of $\Zof\subset\X$. The assertion (1) does not hold by replacing $Z_f$ with its product with a torus. This shows that some {\it ampleness\1} is needed.
\msn
{\bf 2.5.~Proof of the equalities in $(4)$.} By Remark~(2.3), we get the equality in the Grothendieck group of complex algebraic varieties
$$[\Zf]=\mopl_{\si\les\De}\,(-1)^{c_{\si}}[\Zofs],
\leqno(2.5.1)$$
where $\Zofs:=\Zof\cap\Xo_{\si}$, which is the closure of $\Zfs$ in $\Zof$.
(Note that $[\Zofs]=\msum_{\tau\les\si}\,[Z_{\!f,\tau}]$.)
This implies the equality in the Grothendieck group of mixed Hodge structures
$$\chi_c(\Zf)=\mopl_{\si\les\De}\,(-1)^{c_{\si}}\chi_c(\Zofs).
\leqno(2.5.2)$$
Combined with Corollary~2, this implies that
$$h^{n-1,p,0}_c(\Zf)=\msum_{d_{\si}=p+1}\,h^{p,p,0}_c(\Zofs)\q(p\in[1,n{-}1]).
\leqno(2.5.3)$$
The first equality of (4) then follows from (2.4.2).
\sk
For the second equality of (4), we have the spectral sequences in $\MHS$
$$\aligned{}_{\X}E_1^{p,q}=\mopl_{c_{\si}=p}\,H^q(\Xo_{\si})&\Longrightarrow H^{p+q}_c(X_{\De}),\\{}_ZE_1^{p,q}=\mopl_{c_{\si}=p}\,H^q(\Zofs)&\Longrightarrow H^{p+q}_c(\Zf),\endaligned
\leqno(2.5.4)$$
where the direct sums are taken over the faces $\si\les\De$ with $c_{\si}=p$, and moreover the differentials
$${}_{\X}d_1^{\1 p,q}:{}_{\X}E_1^{p,q}\to{}_{\X}E_1^{p+1,q},\q{}_Zd_1^{\1 p,q}:{}_ZE_1^{p,q}\to{}_ZE_1^{p+1,q}$$
are induced by those of the complexes in (2.3.5), see Remark~(2.5) below.
\sk
By (the proof of) Corollary~2, we get
$$\Gr^W_0H^q(\Xo_{\si})=\Gr^W_0H^q(\Zofs)=0\q(\forall\,q>0).
\leqno(2.5.5)$$
(Indeed, the argument for $H^q(\Xo_{\si})$ is similar.)
Since the spectral sequence is defined in $\MHS$, this implies the partial $E_2$-degenerations ($\forall\,p\in\Z$)\,:
$${}_{\X}E_2^{p,0}=\Gr^W_0H_c^p(X_{\De}),\q{}_ZE_2^{p,0}=\Gr^W_0H_c^p(\Zf).
\leqno(2.5.6)$$
\sk
Set
$$\aligned{}_{\X}C^{\ssb}&:={}_{\X}E_1^{\1*,0}=\mopl_{c_{\si}=\ssb}\,H^0(\Xo_{\si}),\\{}_ZC^{\ssb}&:={}_ZE_1^{\1*,0}=\mopl_{c_{\si}=\ssb}\,H^0(\Zofs).\endaligned$$
By the weak Lefschetz type theorem, there are canonical isomorphisms
$${}_{\X}C^p\simto{}_ZC^p\q(\forall\,p<n{-}1),
\leqno(2.5.7)$$
together with the injectivity of the morphism in the following definition
$$A:={\rm Coker}({}_{\X}C^{\1 n-1}\into{}_ZC^{\1 n-1}).$$
\sk
We have
$$\dim_{\Q}H^0(\Zofs)=l^*(\si)\pl 1\q\h{if}\q d_{\si}=1.$$
This implies that
$$\dim_{\Q} A=\mopl_{d_{\si}=1}\,l^*(\si).
\leqno(2.5.8)$$
\sk
Since $X_{\De}\,(=\TT)$ and $\Zf$ are smooth affine varieties of dimension $n$ and $n{-}1$ respectively, we have by (2.5.6) the acyclicity
$$H^p{}_{\X}C^{\ssb}=0\q(p\ne n),\q H^p{}_ZC^{\ssb}=0\q(p\ne n{-}1),
\leqno(2.5.9)$$
together with the equality
$$H^n{}_{\X}C^{\ssb}\cong\Q.
\leqno(2.5.10)$$
\sk
Set
$${}_{\X}B^p:={\rm Im}\,{}_{\X}d^{\1 p-1,0}\subset{}_{\X}C^p,\q{}_ZB^p:={\rm Im}\,{}_Zd^{\1p-1,0}\subset{}_ZC^p.$$
Using (2.5.7), (2.5.9) together with the line lemma, we can get the commutative diagram of short exact sequences
$$\begin{array}{cccccccccccccc}
&& 0 && 0 && 0\\
&&\downarrow&&\downarrow&&\downarrow\\
0&\to&{}_{\X}B^{n-1}&\to&{}_{\X}C^{\1 n-1}&\to & {}_{\X}B^n &\to & 0\\
&&|\!\1 |&&\downarrow&&\downarrow\\
0&\to&{}_ZB^{n-1}&\to &{}_ZC^{\1n-1}&\to & {}_ZE_2^{\1 n-1,0} &\to & 0\\
&&\downarrow&&\downarrow&&\downarrow\\
&& 0 &\to & A\, &=& A\, &\to & 0\\
&&&&\downarrow&&\downarrow\\
&& && 0 && 0
\end{array}$$
Note that the vertical morphism in the left column is an isomorphism by (2.5.7), (2.5.9).
\sk
The equality (2.5.10) gives the short exact sequence
$$0\to{}_{\X}B\1^n\to{}_{\X}C\1^n\to\Q\to 0.$$
This implies the last equality of the following.
$$\#\{\si\les\De\mid d_{\si}=0\}=\dim_{\Q}{}_{\X}C\1^n=\dim_{\Q}
{}_{\X}B\1^n+1.$$
By (2.5.8) we get moreover the equality
$$\Pi=\dim_{\Q}{}_{\X}C\1^n+\dim_{\Q}A.$$
The second equality of (4) now follows from the exactness of the right column of the diagram.
\msn
{\bf Remark~2.5.} Using the {\it classical\1} $t$-structure on $D^b\MHM(\X)$ (see \cite[4.6,2]{mhm}) together with Corollary~(2.3), we can construct filtered bounded complexes of mixed Hodge modules $(\M_{\X}^{(q)\1\ssb},G)$ by increasing induction on $q\in\N$ so that there are isomorphisms in $D^b\MHM(\X)$\,:
$$\aligned&\Gr_G^p\1\M_{\X}^{(q)\1\ssb}\cong\begin{cases}\K_{\X}^{(p)\1\ssb}[-p]&(p\in[0,q]),\\ \,\,0&(p>q).\end{cases}\\&\h{with}\q\K_{\X}^{(p)\1\ssb}:=\mopl_{c_{\si}=p}\,\Q_{h,\Xo_{\si}},\endaligned$$
in the notation of (1.8). Indeed, $\M_{\X}^{(q)\1\ssb}$ is defined by the shifted mapping cone as follows:
$$\M_{\X}^{(q)\1\ssb}:=C\bl(\M_{\X}^{(q-1)\1\ssb}\to\K_{\X}^{(q)\1\ssb}[1{-}q]\br)[-1].$$
The morphism in the mapping cone is the composition of the morphisms of $D^b\MHM(\X)$
$$\M_{\X}^{(q-1)\1\ssb}\to\K_{\X}^{(q-1)\1\ssb}[1{-}q]\to\K_{\X}^{(q)\1\ssb}[1{-}q],$$
where the last morphism induced by the differential of $K_{\X}^{\ssb}$ in Corollary~(2.3) using the adjunction relation for closed immersions $i:Z'\subset Z$\,:
$${\rm Hom}(\Q_{h,Z},i_*\Q_{h,Z'})={\rm Hom}(i^*\Q_{h,Z},\Q_{h,Z'})={\rm End}(\Q_{h,Z'}).$$
The representative of $\K_{\X}^{(q)\1\ssb}$ in $C^b\MHM(\X)$ is replaced in order to define the mapping cone. The filtration $G$ on the mapping cone is defined naturally.
\sk
Using $(\M_{\X}^{(q)\1\ssb},G)$ for $q$ sufficiently large, we get the first spectral sequence for $X_{\De}$ in (2.5.4). The argument is similar for the spectral sequence for $\Zf$, see also \cite{BBD}.
\msn
{\bf 2.6.~Algorithm for the non-prime case.} It does not seem easy to understand the algorithm of Danilov--{Khovanski\u\i} \cite[5.2]{DK} on the Hodge numbers with compact supports of a non-degenerate toric hypersurface $Z_f$ in the non-prime case. It is stated there as follows:
``By 3.11, the numbers $e^{p,q}(Z)$ for $p\1{+}\1 q\1{>}\1 n\mi 1$ are also known; hence the same is true for the numbers $e^{p,q}(\Zo)$ for $p\1{+}\1 q\1{>}\1 n\mi 1$."
Here $\Zo$ is our $\Zo_{\De'}$, and 3.11 seems to mean a weak Lefschetz type theorem for {\it ample\1} divisors.
The meaning of ``hence the same is true" is rather non-trivial. It turns out that we have to use the equality (2.6.2) below (which is a {\it refinement\1} of an equality written before the above sentence in terms of the cohomology with compact supports of toric hypersurfaces indexed by the faces of $\De'$) instead of 3.11, see Remark~(2.6)\,(ii) below.
It is not difficult to see that the compactification $\Zo_{\De'}$ of $\Zf$ in a toric desingularization $\X_{\De'}$ of $\X_{\De}$ may have non-middle cohomology with Hodge level more than 0 using Example~(1.7), see Remark~(2.6)\,(i) below. This shows that the weak Lefschetz type theorem does not hold, and $\Zo_{\De'}$ is not an ample divisor on $\X_{\De'}$.
\msn
{\bf Remarks~2.6}\,(i). Consider for instance the case where $f$ be a sufficiently general deformation of
$$g:=g_1g_2\q\h{with}\q g_1:=(x{+}1)(y{+}1)\pl z,\,\,\,g_2:=u^a\pl v^a\pl 1,$$
with (compact) Newton polyhedron $\Ga(g)$ {\it unchanged\1} ($a\ges 3$). Here $\Ga(g)$ is the convex hull of ${\rm Supp\,g}\subset\R^5$, and $\De=\De_1{\times}\De_2$ with $\De_j=\Ga(g_j)$ ($j=1,2$). Note that $\De_1$ is a special case of $\De$ in Example~(1.7) with $k=4$, and we have the majorizing prime polyhedron $\De'_1$ as in the proof of Theorem~1. This gives the majorizing prime polyhedron $\De':=\De'_1{\times}\De_2$ of $\De$.
\sk
Let $\X_{\De},\X_{\De'}$ be the toric varieties associated with the dual fans of $\De,\De'$ respectively. Let $\Zo_{\De},\Zo_{\De'}$ be the closure of $Z:=g^{-1}(0)\subset(\C^*)^5$ respectively in $\X_{\De},\X_{\De'}$.
We have the proper morphism $\rho:\Zo_{\De'}\to\Zo_{\De}$ together with the decomposition $\X_{\De}=\X_1{\times}\X_2$ associated with $\De=\De_1{\times}\De_2$. Let $\Zo_2\subset\X_2\,(=\PP^2)$ be the hypersurface associated with $g_2$.
By the decomposition theorem together with the calculation in Example~(1.7), we get the following non-canonical isomorphism in $D^b_c(\X_{\De}^{\rm an},\Q)$\,:
$$\RR\rho_*\Q_{\Zo_{\De'}}[5]\cong{\rm IC}_{\Zo_{\De}}\Q\oplus{\rm IC}_{\Zo_2}\Q(-2)[-1]\oplus{\rm IC}_{\Zo_2}\Q(-1)[1].
\leqno(2.6.1)$$
Here $\Zo_2$ is identified with $\{P_1\}{\times}\Zo_2\subset\X$ with $P_1\in\X_1$ the unique non-quasi-smooth point corresponding to the vertex $P_1\in\De_1$ contained in 4 faces of $\De_1$ with dimension 2 as in Example~(1.7). Note that ${\rm IC}_{\Zo_2}\Q=\Q_{\Zo_2}[1]$, since $\Zo_2$ is a {\it smooth curve.}
\sk
The non-canonical isomorphism (2.6.1) can be lifted to the bounded derived category of mixed Hodge modules (since the decomposition theorem holds also in this category). We then see that $H^k(\Zo_{\De'})$ cannot have type $(p,p)$ when $k=5\,{\pm}\,1$.
\msn
{\bf Remarks~2.6}\,(ii). Let $\De'$ be a prime cutting of the Newton polyhedron $\De$ of a non-degenerate Laurent polynomial $f$, see (1.3). We have the equality in the Grothendieck ring of mixed $\R$-Hodge structures $K_0(\MHS_{\1\R})$ (where the $\Q$-structure is forgotten)\,:
$$\chi_c(\Zo_{\De'})=\msum_{\si\les\De}\,\bl(\msum_{\tau\les\De',\1\pi(\tau)=\si}\,(\Lb{-}1)^{d_{\tau}-d_{\si}}\bl)\1\chi_c(Z_{\si}),
\leqno(2.6.2)$$
using the K\"unneth formula for cohomology with compact supports, as is well-known. Some calculation about the {\it pull-back\1} of $f$ is also needed before applying the K\"unneth formula, see for instance \cite[Remark~2.1\,(iv)]{des}. Here $Z_{\si}:=\Zo_{\De}\cap X_{\si}$, and $\Lb$ denotes the class of $\Q(-1)$. This shows that $\chi_c(\Zf)=\chi_c(Z_{\De})$ can be determined by $\chi_c(\Zo_{\De'})$ modulo the $\chi_c(Z_{\si})$ for $\si<\De$ (with $d_{\si}<n)$. Moreover, one can determine $\chi_c(\Zo_{\De'})_{>n-1}$ using (2.6.2) and the inductive argument together with the {\it weak Lefschetz\1} type theorem on $\X_{\De}$ and also the {\it ampleness\1} and {\it stratified smoothness\1} of $\Zo_{\De}\subset\X_{\De}$ as is explained in Remark~(2.4). (We denote by $\theta_{>k}$ the weight$\,{>}\,k$ part of $\theta\in K_0(\MHS_{\1\R})$ in general, which is induced by the exact functor of mixed Hodge structures $H\mapsto H/W_kH$, and similarly for $\theta_{<k}$ using the exact functor $H\mapsto W_{k-1}H$.) This is a variant of an argument informed by A.~Stapledon interpreting \cite{DK}. The formula (2.6.2) does not seem to be mentioned explicitly in {\it loc.\,cit.,} although its underlying formula in terms of the cohomology with compact supports of toric hypersurfaces indexed by the faces of $\De'$ is noted, see also the beginning of this subsection. Since $\Zo_{\De'}$ is quasi-smooth and compact, we can get also $\chi_c(\Zo_{\De'})_{<n-1}$ by duality.
\sk
On the other hand, forgetting the weight filtration $W$, the equality (2.6.2) shows that $\chi'_c(\Zo_{\De'})$ can be determined by $\chi'_c(\Zf)$ modulo the $\chi'_c(Z_{\si})$ for $\si<\De$. Here $\chi'_c$ means that the weight filtration $W$ is forgotten and the Grothendieck ring of filtered vector spaces is used. (This is equivalent to considering the subtotal Hodge numbers in the introduction.) Since $\chi'_c(\Zf)$ is calculated in \cite[4.4]{DK}, we can get $\chi'_c(\Zo_{\De'})$ using (2.6.2) together with the inductive argument. Moreover $\chi_c(\Zo_{\De'})_{>n-1}$ and $\chi_c(\Zo_{\De'})_{<n-1}$ are already known by the above argument. We thus get $\chi_c(\Zo_{\De'})$ and hence $\chi_c(\Zf)$ in $K_0(\MHS_{\1\R})$. (Note that an element of $K_0(\MHS_{\1\R})$ is uniquely determined by its Hodge numbers, classifying simple $\R$-Hodge structures.) Here it is also possible to use $\chi_c(\Zf)_{<n-1}$ instead of $\chi'_c(\Zo_{\De'})$ as in \cite[5.2]{DK}.

\end{document}